\documentstyle[12pt]{amsart}
  \def\to{\longrightarrow}
  \def\lra{\to}

  \def\phi{\varphi}
  \def\epsilon{\varepsilon}
  \def\tilde{\widetilde}

  \def\({\left(}
  \def\){\right)}
  \def\O{{\cal O}}
  
  \def\calo{\O}
  
  \newtheorem{theorem}{Theorem}[section]
  \newtheorem{definition}[theorem]{Definition}
  \newtheorem{remark}[theorem]{Remark}
  \newtheorem{lemma}[theorem]{Lemma}
  
  \newtheorem{proposition}[theorem]{Proposition}

  \newenvironment{varthm*}[1]{\trivlist\item[]\it{\bf #1.}}{\endtrivlist}
  \def\bbP{{\mathbb P}}
  \def\mrom#1{{\rm #1}}
  \def\mult{\mathop{\mrom{mult}}\nolimits}
  \def\Pic{\mathop{\mrom{Pic}}\nolimits}
  
\begin{document}
  \thanks{ Both authors were partially supported by KBN grant 2 P03A 083 10.
    They are fellows of the Foundation for Polish Science.}
  \title{Double covers of ${\mathbb P}^{3}$ and
    Calabi-Yau varieties} 
   \author{S\L awomir\ Cynk} \address{Instytut Matematyki,
    Uniwersytet Jagiello\'nski,\newline
    Reymonta 4, PL-30-059 Krak\'ow,
    Poland}
    \email{cynk@@im.uj.edu.pl\\szemberg@@im.uj.edu.pl} 
    \author{Tomasz\ Szemberg} 
    \begin{abstract}
   We study a class of Calabi--Yau varieties that can be represented as a
   non--singular model of a double covering of ${\mathbb P}^{3}$ branched along
   certain octic surfaces. We compute Euler numbers of all constructed
   examples and describe  their resolution of singularities. 
 \end{abstract}
  \maketitle
  \section*{Introduction.}
  In the present paper we study examples of double coverings of the projective
  space $\bbP^3$ branched over an octic surface.  A double covering of $\bbP^3$
  branched over a smooth octic is a Calabi-Yau threefold. If the octic is
  singular then so is the double covering and we study its resolution of
  singularities.  In this paper we restrict our considerations to the case of
  octics with only non-isolated singularities of a special type, namely looking
  locally like plane arrangements.

  Our research was inspired by a paper of Persson \cite{Per85} where K3
  surfaces arising as double covers of $\bbP^2$ branched over curves of degree
  six are studied. In this note we also adopt some methods introduced in
  \cite{Hun89} by Hunt in studying Fermat covers of $\bbP^3$ branched over
  plane arrangements.

  The main results of this note are Theorem~\ref{blowups} and
  Theorem~\ref{main} which can be formulated together as follows
  \begin{varthm*}{Theorem}
    Let $S\subset\bbP^3$ be an octic arrangement with no $q-$fold curve for
    $q\geq 4$ and no $p-$fold point for $p\geq 6$. Then the double covering of
    $\bbP^3$ branched along $S$ has a non-singular model $Y$ which is a
    Calabi-Yau threefold.

    Moreover if $S$ contains no triple elliptic curves and $l_3$ triple lines
    then the Euler characteristic $e(Y)$ of $Y$ is given as follows \samepage
    $$e(Y)=8-\sum_i(d_i^3-4d_i^2+6d_i)+2\sum_{i\neq j}(4-d_i-d_j)d_id_j-$$
    $$-\sum_{i\neq j\neq k\neq
      i}d_id_jd_k+4p_4^0+3p_4^1+16p_5^0+18p_5^1+20p_5^2+l_3,$$
    where $d_i$
    denotes the degree of the arrangement surfaces and $p_i^j$ the number of
    $i-$fold points contained in $j$ triple curves.
  \end{varthm*}
  For the arrangements we allow exactly six types of singularities.  For each
  case we describe precisely the resolution of singularities in the double
  cover.  Then we study the effect on the Euler number of every blowing up.
  This leads to the formula on the Euler number of $Y$ as stated in the above
  Theorem. Using this formula we obtain a table of examples of Calabi-Yau
  threefolds with $63$ different Euler numbers. In the view of the mirror
  symmetry it is important to have examples of Calabi-Yau threefolds with as
  many Euler numbers as possibly (it is conjectured that there are only
  finitely many possible numbers to appear as the Euler number of a Calabi-Yau
  threefold).
    \section{Admissible blowing ups.}
    Let $D$ be a reduced divisor on a smooth threefold $X$. We assume that $D$
    is even as an element of the Picard group $\Pic X$, this means that there
    exists $B\in \Pic X$ such that $D\cong 2B$. In this case there exists a
    double covering $\pi :Y\lra X$ branched along $D$. If the divisor is smooth
    then $Y$ is also smooth and we have
    \begin{itemize}
    \item[] $h^i(\calo_Y)=h^i(\calo_X)+h^i(\calo_X(-B)),$
    \item[] $K_Y\cong\pi^*(K_X+B)$,
    \item[] $e(Y)=2e(X)-e(D)$.
    \end{itemize}
    If $D$ is singular then the type of its singularities determines
    singularities of $Y$. Hence it is enough to consider an embedded resolution
    $D$' of $D$ to obtain a resolution for $Y$. The problem however is to
    ensure that the resulting divisor is even.

    Persson, studying double sextics, introduced in \cite{Per85} a notion of
    inessential singularities i.e.\ such which do not affect the Euler
    characteristic and the canonical divisor of the double cover. In a
    threedimensional case an analogous characterization would be much more
    complicated. Moreover in our paper we are interested in Calabi-Yau
    manifolds with many different Euler numbers.  This leads to the following
    description of blowing-ups which do not affect the first Betti number and
    the canonical divisor of the double cover.

    Let $X$ be a smooth threefold and $D\subset X$ an even, reduced divisor.
    Let $Z\subset D$ be a smooth irreducible proper subvariety and let $\sigma
    :\tilde{X}\lra X$ be the blowing up of $X$ in $Z$ with an exceptional
    divisor $E$. By $\mult_{Z/D}$ we denote the generic multiplicity of $D$ at
    $Z$ and by $\tilde{D}$ the proper transform of $D$.  Then the divisor
    $D^*\subset\tilde{X}$ defined as
    $$D^*:=\left\{\begin{tabular}{ll}
        $\tilde{D}$ & if  $\mult_{Z/D}$ is even, \\
        $\tilde{D}+E$ & if $\mult_{Z/D}$ is odd
                  \end{tabular}\right.$$
                is the only reduced and even divisor satisfying $\tilde{D}\leq
                D^*\leq\sigma^*D$.
  \begin{definition}
    Let $Z\subset D\subset X$ be as above. We call the blowing up
    $\sigma:\tilde{X}\lra X$ admissible iff
    $$K_{\tilde{X}}+\frac{1}{2}D^*\cong\sigma^*(K_X+\frac{1}{2}D).$$
  \end{definition}
  The following proposition gives a characterization of admissible blowing ups.
  \begin{proposition}
    On a smooth threefold $X$ there are exactly four types of admissible
    blow-ups:
    \begin{itemize}
      \leftskip=3cm
    \item[C)] blowing up of a curve $Z$ with $\mult_{Z/D}=2 \mbox{ or } 3$,
    \item[P)] blowing up a point $Z$ with $\mult_{Z/D}=4 \mbox{ or } 5$.
    \end{itemize}
  \end{proposition}
  \proof Let $r$ denote the codimension of $Z$ in $X$ and let $m=mult_{Z/D}$.
  Furthermore let $\varepsilon$ be equal $0$ if $m$ is even and $1$ otherwise.
  Then we have $K_{\tilde{X}}\cong\sigma^*K_X+(r-1)E$ and
  $D^*=\sigma^*D-(m-\varepsilon)E$. Therefore

  $$K_{\tilde{X}}+\frac{1}{2}D^*\cong\sigma^*(K_X+D)+(r-1-\frac{m-\varepsilon}{2}))E.$$
  It follows that $\sigma$ is admissible if and only if $m=2(r-1)+\varepsilon$.
  As the only solutions for $r=2$ we get $m=2$ or $m=3$ and for $r=3$
  respectively $m=4$ or $m=5$.  \qed
  \section{Resolution of singularities of the branch locus.}
  Before we start to resolve singularities we need some definitions.  We
  consider octic surfaces in $S\subset\bbP^3$ which locally look like an
  arrangement of planes (see \cite[p.109]{Hun89}).  More precisely let $S$ be a
  sum of smooth irreducible surfaces $S_1,\ldots,S_r$ contained in a smooth
  threefold $U$.  We suppose that
    \begin{itemize}
    \item[1.] For any $i\neq j$ the surfaces $S_i$ and $S_j$ intersect
      transversally along a smooth irreducible curve $C_{i,j}$ or they are
      disjoint,
    \item[2.] The curves $C_{i,j}, C_{k,l}$ either coincide or they intersect
      transversally,
    \end{itemize}
    We call a surface satisfying above conditions an {\it arrangement}.  If
    $U=\bbP^3$ and surfaces $S_1,\ldots,S_r$ are of degree $d_1,\ldots,d_r$
    respectively, with $d_1+\cdots+d_r=8$ then we call $S$ an {\it octic
      arrangement}.

    We say that an irreducible curve $C\subset S$ is a $q-$fold curve if
    exactly $q$ of surfaces $S_1,\ldots,S_r$ pass through it. A point $P\in S$
    is called a $p-$fold point if exactly $p$ of the surfaces $S_1,\ldots,S_r$
    pass through it.
  \begin{theorem}\label{blowups}

    Let $S\subset\bbP^3$ be an octic arrangement with no $q-$fold curve for
    $q\geq 4$ and no $p-$fold point for $p\geq 6$. Then there exists a sequence
    $\sigma=\sigma_1\circ\ldots\circ\sigma_s:\bbP^*_3\lra\bbP^3$ of admissible
    blowing ups and a smooth and even divisor $S^*\subset\bbP^*_3$ such that
    $\sigma_*(S^*)=S$ and the double covering $Y$ of $\bbP^*_3$ branched over
    $S^*$ is a smooth Calabi-Yau threefold.
  \end{theorem}  
  \proof
  The proof consists of careful resolving the singularities of $S$.\\
  {\bf a) Fivefold points.}  First we blow up $\sigma_1:\bbP^3_{(1)}\lra\bbP^3$
  all points $P\in\bbP^3$ which are $5-$fold points of $S$.  Then in the
  exceptional divisor $\bbP^2$ over each of these points we get one of the
  following configurations of lines - exceptional divisors for surfaces passing
  through $P$.\\\parbox{4.2cm}{\begin{picture}(140,90)(-20,-20)
     \put(0,0){\line(1,0){90}}
     \put(0,0){\line(0,1){60}}
     \put(90,0){\line(0,1){60}}
     \put(0,60){\line(1,0){90}}
     \put(0,10){\line(2,1){90}}
     \put(0,50){\line(2,-1){90}}
     \put(0,40){\line(1,0){90}}
     \put(60,0){\line(1,3){20}}
     \put(80,0){\line(-1,3){20}}
     \put(20,40){\circle*{1}}
     \put(40,30){\circle*{1}}
     \put(65.7,17){\circle*{1}}
     \put(65.7,43){\circle*{1}}
     \put(60,40){\circle*{1}}
     \put(70,30){\circle*{1}}
     \put(66.7,40){\circle*{1}}
     \put(73.3,40){\circle*{1}}
     \put(76,48){\circle*{1}}
     \put(76,12){\circle*{1}}
     \put(25,46){\small$\bbP^2$}
     \put(0,-10){\small Case a-1}
     \end{picture}}\parbox{4.2cm}{
     \begin{picture}(140,90)(-20,-20)
     \put(0,0){\line(1,0){90}}
     \put(0,0){\line(0,1){60}}
     \put(90,0){\line(0,1){60}}
     \put(0,60){\line(1,0){90}}
     \put(0,10){\line(2,1){90}}
     \put(0,50){\line(2,-1){90}}
     \put(0,20){\line(4,1){90}}
     \put(60,0){\line(1,3){20}}
     \put(80,0){\line(-1,3){20}}
     \put(72.7,38.2){\circle*{1}}
     \put(40,30){\circle*{1.5}}
     \put(65.7,17){\circle*{1}}
     \put(65.7,43){\circle*{1}}
     \put(67.7,36.9){\circle*{1}}
     \put(70,30){\circle*{1}}
     \put(76,48){\circle*{1}}
     \put(76,12){\circle*{1}}
     \put(25,46){\small$\bbP^2$}
     \put(0,-10){\small Case a-2}    
     \end{picture}}\parbox{4.2cm}{
     \begin{picture}(140,90)(-20,-20)
     \put(0,0){\line(1,0){90}}
     \put(0,0){\line(0,1){60}}
     \put(90,0){\line(0,1){60}}
     \put(0,60){\line(1,0){90}}
     \put(0,10){\line(2,1){90}}
     \put(0,50){\line(2,-1){90}}
     \put(0,30){\line(1,0){90}}
     \put(60,0){\line(1,3){20}}
     \put(80,0){\line(-1,3){20}}
     \put(40,30){\circle*{2}}
     \put(40,30){\circle*{1}}
     \put(65.7,17){\circle*{1}}
     \put(65.7,43){\circle*{1}}
     \put(70,30){\circle*{1}}
     \put(76,48){\circle*{1}}
     \put(76,12){\circle*{1}}
     \put(25,46){\small$\bbP^2$}
     \put(0,-10){\small Case a-3}        
     \end{picture}
     }

   Obtained cases depend on the number of threefold curves passing through $P$.
   Now we replace $S$ by its strict transform plus the exceptional divisor. And
   we call this new branch locus $S_1$. If there were no $5-$fold points on $S$
   we have
   $S_1=S$. We observe that $S_1$ contains no $5-$fold points.\\
   \parbox{8cm}{ {\bf b) Triple curves.}  If there are triple curves on $S_1$
     we blow them up $\sigma_2:\bbP^3_{(2)}\lra\bbP^3_{(1)}$. Let $C\subset
     S_1$ be a triple curve. Then in the exceptional divisor $C\times\bbP^1$ we
     get the following configuration where $C_1, C_2, C_3$ are isomorphic to
     $C$ and $L_1,\ldots ,L_t$ are lines,
   with $t$ equal to the number of $4-$fold points on $C$.}%
 \parbox{4.7cm}
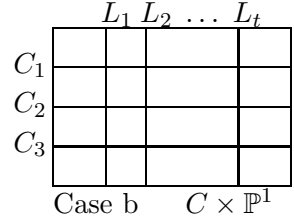
{%
     \begin{picture}(140,90)(-30,-20)
     \put(0,0){\line(1,0){90}}
     \put(0,0){\line(0,1){60}}
     \put(90,0){\line(0,1){60}}
     \put(0,60){\line(1,0){90}}
     \put(0,15){\line(1,0){90}}
     \put(0,30){\line(1,0){90}}
     \put(0,45){\line(1,0){90}}
     \put(20,0){\line(0,1){60}}
     \put(35,0){\line(0,1){60}}
     \put(70,0){\line(0,1){60}}
     \put(-15,13){\small$C_{3}$}
     \put(-15,28){\small$C_{2}$}
     \put(-15,43){\small$C_{1}$}
     \put(18,62){\small$L_{1}$}
     \put(33,62){\small$L_{2}$}
     \put(50,62){\small$\dots$}
     \put(68,62){\small$L_{t}$}
     \put(20,15){\circle*{1}}
     \put(35,15){\circle*{1}}
     \put(70,15){\circle*{1}}
     \put(20,30){\circle*{1}}
     \put(35,30){\circle*{1}}
     \put(70,30){\circle*{1}}
     \put(20,45){\circle*{1}}
     \put(35,45){\circle*{1}}
     \put(70,45){\circle*{1}}
     \put(50,-10){\small$C\times\bbP^1$}
     \put(0,-10){\small Case b}    
     \end{picture}
     }\\
   As a new branch locus for the double cover we take $S_2$ equal to the strict
   transform of $S_1$ plus the exceptional divisors. There are no $5-$fold
   points
   and no triple curves on $S_2$.\\
   \parbox{8cm}{ {\bf c) Fourfold points.}  Now we blow up
     $\sigma_3:\bbP^3_{(3)}\lra\bbP^3_{(2)}$ all the $4-$fold points of $S_2$.
     In the exceptional divisor over a $4-$fold
   point $P\in S_2$ we get the following configuration.}%
 \parbox{4.7cm}{%
     \begin{picture}(140,90)(-30,-20)
     \put(0,0){\line(1,0){90}}
     \put(0,0){\line(0,1){60}}
     \put(90,0){\line(0,1){60}}
     \put(0,60){\line(1,0){90}}
     \put(0,0){\line(2,1){90}}
     \put(14,0){\line(3,5){36}}
     \put(0,40){\line(2,-1){80}}
     \put(35,60){\line(1,-2){30}}
     \put(20,10){\circle*{1}}
     \put(60,10){\circle*{1}}
     \put(40,20){\circle*{1}}
     \put(41.8,46.4){\circle*{1}}
     \put(52,26){\circle*{1}}
     \put(29.2,25.4){\circle*{1}}
     \put(15,46){\small$\bbP^2$}
     \put(0,-10){\small Case c}
       \end{picture}
 }\\
 Now as $S_3$ we take the strict transform of $S_2$.
 On $S_3$ there are no more singularities of types a), b), c).\\
 \parbox{8cm}{ {\bf d) Double curves.}  In the last step we blow up
   $\sigma_4:\bbP_3^*=\bbP^3_{(4)}\lra\bbP^3_{(3)}$ the double curves on $S_3$
   and take $S^*=S_4$ as the strict transform of $S_3$.  Since $S$ was an octic
   arrangement we get $S^*$ smooth and even as an element
   of $\Pic\bbP^*_3$.}%
 \parbox{4.7cm}{%
    \begin{picture}(140,90)(-30,-20)
    \put(0,0){\line(1,0){90}}
    \put(0,0){\line(0,1){60}}
    \put(90,0){\line(0,1){60}}
    \put(0,60){\line(1,0){90}}
    \put(0,15){\line(1,0){90}}
    \put(0,45){\line(1,0){90}}
    \put(20,0){\line(0,1){60}}
    \put(35,0){\line(0,1){60}}
    \put(70,0){\line(0,1){60}}
    \put(-15,13){$C_{2}$}
    \put(-15,43){$C_{1}$}
    \put(18,62){$L_{1}$}
    \put(33,62){$L_{2}$}
    \put(50,62){$\dots$}
    \put(68,62){$L_{t}$}
    \put(20,15){\circle*{1}}
    \put(35,15){\circle*{1}}
    \put(70,15){\circle*{1}}
    \put(20,45){\circle*{1}}
    \put(35,45){\circle*{1}}
    \put(70,45){\circle*{1}}
    \put(0,-10){\small Case d}
    \end{picture}
}  \\
  
Let $\sigma=\sigma_4\circ\ldots\circ\sigma_1$ and let $\pi:Y\lra\bbP^*_3$ be
the double covering branched over $S^*$. Using the adjunction formula and the
Serre duality we get
$$K_Y\cong\pi^*(K_{\bbP^*_3}+\frac{1}{2}S^*)\cong\pi^*(\sigma^*(K_{\bbP^3}+
\frac{1}{2}S))=\calo_Y$$
and
$$h^1(\calo_Y)=h^1(\calo_{\bbP^*_3})+h^1(\calo_{\bbP^*_3}(-\frac{1}{2}S^*))=$$
  \[\rule{2cm}{0cm}h^1(\calo_{\bbP^*_3})+h^2(\calo_{\bbP^*_3}(K_{\bbP^*_3}+\frac{1}{2}S^*))=
    h^1(\calo_{\bbP^3})+h^2(\calo_{\bbP^3})=0.%  
    \raisebox{-.3cm}{\rule{2cm}{0cm}\qedsymbol} \]
\begin{remark}
  Observe that the first three steps of the above resolution are uniquely
  determined whereas the last step is defined only upto the order in which we
  blow-up the double curves. Change in the order in which the double curves
  are blown-up may lead to a flop of the resulting threefold, see e.g.
  \cite{kol}. However this does not affect the Euler number in which we are
  interested. Hence we do not distinguish birational models differing by a
  flop.
\end{remark}
\section{Euler characteristic of double octics.}
In this section we compute the Euler characteristic of Calabi-Yau threefolds
obtained from octic arrangements as in Theorem~\ref{blowups}.

For an arrangement $S$ we introduce the following notation:

\addvspace{2em}

  \begin{tabular}{ll}
    $e^*(S)$& is the sum of Euler numbers of all components of $S$,\\ &\\
    $E_p(S)$& is the sum of Euler numbers of $p-$fold curves on $S$,\\ &\\
    $p_j(S)$& is the number of isolated $j-$fold points on $S$,\\ &\\
    $p_j^k(S)$& is the number of isolated $j-$fold points
          lying on exactly $k$ \\& triple curves.\\ &
  \end{tabular}\\
  Now we compute how the above data changes under blowing ups described in the
  proof of Theorem~\ref{blowups}.
\begin{proposition}\label{invariant}
  
  Let $S\subset U$ be an arrangement in a threefold $U$.  Let $\sigma:V\lra U$
  be a blowing up of the type a), b), c) or d) with the center $Z$ and the
  exceptional divisor $E$.  As before $S^*=\sigma^*S+\varepsilon E$, with
  $\varepsilon =0\mbox{ or } 1$ depending on the case a-d).  In this situation
  we have
  $$2e(U)-e^{*}(S)+2E_2(S)-p_3(S)+6E_3(S)+12p^2_5(S)+9p^1_5(S)+6p^0_5(S)=$$
  $$=2e(V)-e^{*}(S^*)+2E_2(S^*)-p_3(S^*)+6E_3(S^*)+12p^2_5(S^*)+9p^1_5
  (S^*)+6p^0_5(S^*)$$
\end{proposition}
\proof The following table describes how the blowing up affects the Euler
numbers of the threefold and the arrangement. The next table shows changes of
the combinatorial data. If $Z$ is a $q-$fold line then $t$ denotes the number
of $(q+1)-$fold points on $Z$.
  \begin{center}
    \begin{tabular}{|c|c|c|c|c|}
    \hline
    type& $e(V)-e(U)$& $e^{*}(S^*)-e^{*}(S)$& $E_2(S^*)-E_2(S)$& 
    $E_3(S^*)-E_3(S)$\\
    \hline
    a-1& $2$& $8$& $10$ &$0$\\
    \hline
    a-2& $2$& $8$& $10$& $0$\\
    \hline
    a-3& $2$& $8$& $10$& $0$\\
    \hline
    b& $e(Z)$& $2e(Z)+t$& $3e(Z)+2t$& $-e(Z)$\\
    \hline
    c& $2$& $4$& $0$& $0$\\
    \hline
    d& $e(Z)$& $t$& $-e(Z)$& $0$ \\
    \hline
\multicolumn{5}{c}{}\\
    \hline
    type& $p_3(S^*)-p_3(S)$& $p_5^0(S^*)-p_5^0(S)$& $p_5^1(S^*)-p_5^1(S)$&
          $p_5^2(S^*)-p_5^2(S)$\\
    \hline
    a-1& $10$& $-1$& $0$& $0$\\
    \hline
    a-2& $7$& $0$& $-1$& $0$\\
    \hline
    a-3& $4$& $0$& $0$& $-1$\\
    \hline
    b& $3t$& $0$& $0$& $0$\\
    \hline
    c& $0$& $0$& $0$& $0$\\
    \hline
    d& $-t$& $0$& $0$& $0$\\
    \hline
    \end{tabular}
  \end{center}
  It is simple to check the entries of the above tables and then case by case
  to verify that the term given in the proposition remains invariant.  
  \qed\medskip 

  As a corollary from the above Proposition we obtain the following
\begin{proposition}\label{eulerupstairs}
  
  In the setup of Theorem~\ref{blowups} we have
  $$e(Y)=8-e^{*}(S)+2E_2(S)-p^0_3(S)+6E_3(S)+12p^2_5(S)+9p_5^1(S)+6p_5^0(S).$$
\end{proposition}
Now we give some formula to compute the invariants used in the above
Proposition.
\begin{lemma}\label{eulersc}
  
  For an octic arrangement $S$ in $\bbP^3$ we have
  $$e^*(S)=\sum d_i^3-4d_i^2+6d_i, \mbox{ and }$$
  $$2E_2(S)+6E_3(S)=2-\sum_{i\neq j}(4-d_i-d_j)d_id_j.$$
\end{lemma}
\proof From the adjunction formula we have
$$e(X)=d^3-4d^2+6d$$
for a smooth surface $X\subset\bbP^3$ of degree $d$.
Similarly if $C\subset\bbP^3$ is a smooth complete intersection of surfaces of
degree $d_1, d_2$ we have
$$e(C)=(4-d_1-d_2)d_1d_2.$$
\qed In \cite{Hun89} a formula for the number of
singular points of an arrangement of planes is given. This formula can be
generalized to the case of an arbitrary octic arrangement. In the simplest
case if there are no triple curves we have
$$\sum_{q\geq 3}{q\choose 3}p_q=\sum_{i\neq j\neq k\neq i}d_id_jd_k.$$
Situation becomes more complicated if there are also triple curves in the
arrangement. Since $\deg(S)=8$ there are only two possibilities:
  \begin{itemize}
  \item either there is one triple elliptic curve and no more triple curves
  \item or there are only triple lines.
  \end{itemize}
  We can easily classify arrangements with a triple elliptic curve
  \begin{itemize}
  \item $d_1=d_2=d_3=d_4=2, \ \ p_3=8, \ \ p_5=0, \ \ e=-16$;
  \item $d_1=d_2=d_3=2, \ \ d_4=d_5=1, \ \ p_3=6, \ \ p_5=0, \ \ e=12$;
  \item $d_1=d_2=d_3=2, \ \ d_4=d_5=1, \ \ p_3=3, \ \ p_5^0=p_5^2=0, \ \ 
    p_5^1=1, \ \ e=24$;
  \item $d_1=d_2=d_3=2, \ \ d_4=d_5=1, \ \ p_3=0, \ \ p_5^0=p_5^2=0, \ \ 
    e=36$;
  \end{itemize}
\begin{lemma}\label{formula}  
  For an arrangement with $l_3$ triple lines ($l_3=\frac{1}{2}E_3$) and no
  triple elliptic curves we get the following formulas:
  $$p_3+4p_4+10p_5-(p_4^1+p_5^1+2p_5^2-l_3)=\sum_{i\neq j\neq k\neq
    i}d_id_jd_k,$$
  $$5l_3=p_4^1+2p_5^1+4p_5^2.$$
\end{lemma}
\proof The number on the right hand side of the first equation is just the sum
of intersection numbers of all possible triples of arrangement surfaces i.e.
the number of triple points in case all intersections are transversal and
reduced.  On the left hand side we take account of the multiple points and
multiple lines. For example if there is a triple line in the picture then it
is an interesection of three planes (cf. condition 1 for the arrangement) and
corresponds to one "lost" point. A $4-$fold point corresponds to $4$ "lost"
points etc.

For the second formula we consider triple lines. They lie on three planes and
intersect the remaining quintic. On the left hand side there is an expected
number of points and on the right hand side the number with multiplicities.
\qed

From Proposition~\ref{eulerupstairs}, Lemma~\ref{eulersc} and
Lemma~(\ref{formula}) we get the main result of this paper.
\begin{theorem}\label{main}
  
  If $S$ is an octic arrangement with no triple elliptic curve then
  $$e(Y)=8-\sum_i(d_i^3-4d_i^2+6d_i)+2\sum_{i\neq j}(4-d_i-d_j)d_id_j-$$
  $$-\sum_{i\neq j\neq k\neq
    i}d_id_jd_k+4p_4^0+3p_4^1+16p_5^0+18p_5^1+20p_5^2+l_3$$
\end{theorem}
\section{Examples.}
In this section we apply the preceding results to various examples of octic
arrangements and for the double covering we get 63 distinct Euler numbers. We
collect the data describing an arrangement and the corresponding Euler number.
Some Euler numbers arise from more than one arrangement, the table contains
only one example in order to keep the paper reasonable short.  \newpage
  \begin{displaymath}
\begin{array}[t]{|l||c|c|c|c|c|c||c|}
\hline
(d_{1},\dots,d_{r})&{p_{4}}^{0}&{p_{4}}^{1}&{p_{5}}^{0}&{p_{5}}^{1}%
&{p_{5}}^{2}&{l_3}&{\bf e(Y)}\rule[-1ex]{0ex}{3.5ex}\\
\hline\hline
8&&&&&&&-296\\
\hline
1,7&&&&&&&-240\\
\hline
2,6&&&&&&&-200\\
\hline
1,1,6&&&&&&&-180\\
\hline
3,5&&&&&&&-176\\
\hline
4,4&&&&&&&-168\\
\hline
1,2,5&&&&&&&-140\\
\hline
1,1,1,5&&&&&&&-120\\
\cline{2-8}
&1&&&&&&-116\\
\hline
2,2,4&&&&&&&-104\\
\hline
2,3,3&&&&&&&-92\\
\hline
1,1,2,4&&&&&&&-84\\
\cline{2-8}
&1&&&&&&-80\\
\cline{2-8}
&2&&&&&&-76\\
\hline
1,1,3,3&&&&&&&-72\\
\cline{2-8}
&1&&&&&&-68\\
\hline
1,1,1,1,4&&&&&&&-64\\
\cline{2-8}&1&&&&&&-60\\
\hline
1,2,2,3&&&&&&&-56\\
\cline{2-8}&1&&&&&&-52\\
\hline
1,1,1,1,4&&&1&&&&-48\\
\hline
1,2,2,3&3&&&&&&-44\\
\hline
2,2,2,2&&&&&&&-40\\
\hline
1,1,1,2,3&&&&&&&-36\\
\cline{2-8}
&1&&&&&&-32\\
\cline{2-8}&2&&&&&&-28\\
\cline{2-8}&3&&&&&&-24\\
\hline 
1,1,2,2,2&&&&&&&-20\\
\hline
1,1,1,1,1,3&&&&&&&-16\\
\cline{2-8}&1&&&&&&-12\\
\cline{2-8}&2&&&&&&-8\\
\cline{2-8}&3&&&&&&-4\\
\hline 
1,1,1,1,2,2&&&&&&&0\\
\cline{2-8}&1&&&&&&4\\
\cline{2-8}&2&&&&&&8\\
\cline{2-8}&3&&&&&&12\\
\cline{2-8}&&&1&&&&16\\
\hline
\end{array}
\end{displaymath}
\vfil \newpage
\begin{displaymath}
  \begin{array}[t]{|l||c|c|c|c|c|c||c|}
\hline
(d_{1},\dots,d_{r})&{p_{4}}^{0}&{p_{4}}^{1}&{p_{5}}^{0}&{p_{5}}^{1}%
&{p_{5}}^{2}&{l_3}&{\bf e(Y)}\rule[-1ex]{0ex}{3.5ex}\\
\hline

\hline 1,1,1,1,1,1,2&&&&&&&20\\
\cline{2-8}&1&&&&&&24\\
\cline{2-8}&2&&&&&&28\\
\cline{2-8}&3&&&&&&32\\
\cline{2-8}&&&1&&&&36\\
\hline 
1,1,1,1,1,1,1,1 &&&&&&&40\\
\cline{2-8}&1&&&&&&44\\
\cline{2-8}&2&&&&&&48\\
\cline{2-8}&3&&&&&&52\\
\cline{2-8}&4&&&&&&56\\
\cline{2-8}&5&&&&&&60\\
\cline{2-8}&6&&&&&&64\\
\cline{2-8}&7&&&&&&68\\
\cline{2-8}&8&&&&&&72\\
\cline{2-8}&9&&&&&&76\\
\cline{2-8}&&1&&2&&1&80\\
\cline{2-8}&&8&&1&&2&84\\
\cline{2-8}&12&&&&&&88\\
\cline{2-8}&&4&&1&1&2&92\\
\cline{2-8}&&6&&2&&2&96\\
\cline{2-8}&&7&&&2&3&104\\
\cline{2-8}&&9&&1&1&3&108\\
\cline{2-8}&&3&&&3&3&112\\
\cline{2-8}&1&3&&&3&3&116\\
\cline{2-8}&2&3&&&3&3&120\\
\cline{2-8}&&4&&&4&4&136\\
\hline
  \end{array}
\end{displaymath}

Most of arrangements from the above table are easy to construct.  We conclude
this paper by giving hints as to construct those
which, at least to the authors, seem less obvious.\\
{\bf $e=-44$}\\
We take two general quadrics and a plane. They have four common points.  Now,
we add a cubic passing through three of these points. Bertini Theorem implies
that the transversality conditions are satisfied away
of the fourth point. But there they are obviously fulfilled.\\
{\bf $e=-24$}\\
We take two general planes and a cubic. They have three common points.  Then
we take a general quadric through two of these points and a plane
through the third one.\\
{\bf $e=52,56,60,64,68,72,76$}\\
We take six planes as the faces of the cube. Then we let another two
planes pass through $0,1,2,3,4,5,6$ vertices, respectively.\\
{\bf $e=88$}\\
Here we consider faces of the octahedron. They have four $4-$fold
points in the affine part and eight more at the infinity.\\
{\bf $e=104$}\\
We take four planes $P_1, P_2, P_3, P_4$ in the general position.  Let $P_5$
be a general plane through the line $P_1\cap P_2$ and similarly $P_6, P_7$
general planes through the lines $P_2\cap P_3$
and $P_3\cap P_4$. The eighth plane can be taken generally.\\
{\bf $e=112, 116, 120$}\\
We consider two tetrahedra glued by a face and take their seven faces. The two
vertices not lying on the common face
are triple points. Then we add a plane through none, one or two of them.\\
{\bf $e=136$}\\
We take four planes $P_1, P_2, P_3, P_4$ in the general position and add a
general plane through each of the lines $P_1\cap P_2$, $P_1\cap P_3$, $P_2\cap
P_4$, $P_3\cap P_4$.

\medskip

Finally, we remark that many more examples of Euler numbers can be obtained
allowing additionally to the nonisolated singularities some isolated ones,
e.g.\ double points. Similarly one can weaken the assumption on the
reducibility of the arrangement allowing surfaces with selfintersections.
These aspects will be studied elsewhere.  

\end{document}